\documentclass[11pt]{amsart}

\usepackage{a4wide}
\usepackage{amssymb}
\usepackage{amstext}
\usepackage{hyperref}

\let\ge=\geqslant
\let\geq=\geqslant
\let\le=\leqslant
\let\leq=\leqslant

\let\ptl=\partial
\let\Sg=\Sigma
\let\sg=\sigma
\let\eps=\varepsilon
\let\Om=\Omega

\newcommand{\rr}{\mathbb{R}}

\newcommand{\pp}{\mathcal{P}}
\newcommand{\nn}{\mathbb{N}}
\newcommand{\sph}{\mathbb{S}}

\newcommand{\escpr}[1]{\left< #1\right>}
\newcommand{\cone}{{\times\!\!\!\!\times}}
\newcommand{\wt}{\widetilde}

\DeclareMathOperator{\vol}{vol}

\newtheorem{theorem}{Theorem}[section]
\newtheorem{proposition}[theorem]{Proposition}
\newtheorem{lemma}[theorem]{Lemma}
\newtheorem{corollary}[theorem]{Corollary}

\theoremstyle{definition}

\newtheorem{remark}[theorem]{Remark}
\newtheorem{example}{Example}    
\newtheorem{problem}{Problem}

\theoremstyle{remark}

\newenvironment{enum}{\begin{enumerate}
}{\end{enumerate}}

\begin{document}

\title[Isoperimetric regions inside Euclidean cones]{Existence and
characterization of regions \\ minimizing perimeter under a volume
constraint \\
inside Euclidean cones}

\author[M.~Ritor\'e]{Manuel Ritor\'e}
\address{Departamento de Geometr\'{\i}a y Topolog\'{\i}a \\
Universidad de Granada \\ E--18071 Granada \\ Espa\~na}
\email{ritore@ugr.es}

\author[C.~Rosales]{C\'esar Rosales}
\address{Departamento de Geometr\'{\i}a y Topolog\'{\i}a \\
Universidad de Granada \\ E--18071 Granada \\ Espa\~na}
\email{crosales@ugr.es}

\date{February 20, 2003}

\thanks{Both authors have been supported by MCyT-Feder research
project BFM2001-3489}


\begin{abstract}
We study the problem of existence of regions separating a given amount
of volume with the least possible perimeter inside a Euclidean cone. 
Our main result shows that nonexistence for a given volume implies
that the isoperimetric profile of the cone coincides with the one of
the half-space.  This allows us to give some criteria ensuring
existence of isoperimetric regions: for instance, local convexity of
the cone at some boundary point.

We also characterize which are the stable regions in a convex cone,
i.e., second order minima of perimeter under a volume constraint. 
From this it follows that the isoperimetric regions in a convex cone
are the euclidean balls centered at the vertex intersected with the
cone.
\end{abstract}

\maketitle

\thispagestyle{empty}

\section{Introduction}
Let $M\subset\rr^{n+1}$ be a solid cone, a cone over a connected, open
subset $C$ with smooth boundary of the unit sphere $\sph^n$.  In this
work we shall study if it is possible to {\em separate\/} a region
$\Om$ of given volume in $M$ with the least possible perimeter.  This
means that
\[
\mathcal{P}(\Om,M)\le \mathcal{P}(E,M),
\]
for all $E\subset M$ such that $\vol(E)=\vol(\Om)$.  Here
$\mathcal{P}(\,\cdot\,, M)$ denotes the perimeter relative to $M$ and
$\text{vol}(\cdot)$ the $(n+1)$-dimensional Hausdorff measure in
$\rr^{n+1}$.  If such an $\Om$ exists we shall call it an {\em
isoperimetric region\/} of volume $\text{vol}(\Om)$.

We remark that only the boundary area of $\Om$ \emph{inside} $M$
contributes to the perimeter, while $\ptl\Om\cap\ptl M$ does not. 
Hence our problem is quite different from the one of {\em enclosing\/}
a given volume with the least possi\-ble perimeter (the area of the
hypersurface $\ptl\Om\cap\ptl M$ is now taken into account), that has
been recently studied by E.~Stredulinsky and W.~P.~Ziemer~\cite{szi}
and C.~Rosales~\cite{rosales} inside a convex body.  Since a cone $M$
contains round balls of any radius, the isoperimetric inequality in
Euclidean space implies that round balls are the solutions to this
second problem.

The problem of separating a given volume inside a convex body in
Euclidean space has been studied by P.~Sternberg and
K.~Zumbrun~\cite{zumbrun2}, \cite{zumbrun1} and A.~Ros and
E.~Vergasta~\cite{rosvergasta}.  Existence for this problem is
guaranteed by the compactness of the ambient manifold.  In these
papers some properties of the minimizing sets, such as connectedness
of the boundary and estimates on the genus of the free boundary in the
three-dimensional case, are obtained.

To decide whether isoperimetric regions exist or not in a noncompact
manifold is not an easy problem.  The only general results, up to our
knowledge, have been proved by F.~Morgan (see
\cite{morgan}, Theorem~13.4, and the remark below) for homogeneous
manifolds, and by M.~Ritor\'e~\cite{ritore2} for convex surfaces. 
There exist examples \cite{ritore1} of complete surfaces of
revolution for which there are no isoperimetric regions of any given
area.  In \cite{ritore2} it was shown that the direct method of the
calculus of variations (take a minimizing sequence and extract a
convergent subsequence) cannot be applied in the noncompact setting
since part or all of a minimizing sequence could diverge.

We start this work by proving a general result,
Theorem~\ref{th:convergence}, on the behavior of minimizing sequences. 
They can be broken in two parts: passing to a subsequence, one
converges to an isoperimetric region for the volume it encloses, and
the other one diverges.  In some cases it is easier to allow part of
the sequence to diverge in order to get a contradiction, as was done
in \cite{ritore2} for complete convex surfaces.  A survey of known
regularity results completes Section~\ref{sec:existence} of the paper.

In Section~\ref{sec:cone} we treat the existence problem of
isoperimetric regions in a cone $M\subset\rr^{n+1}$.  A simple
observation, due to the presence of dilations, allows us to conclude
that the isoperimetric dimension \cite{gromov}, \S~6.4, of the cone is
$(n+1)$, and that existence of isoperimetric regions for one given
volume implies existence for all volumes.  Since the cone is
asymptotically flat, it is not difficult to show that the
isoperimetric profile of the cone is bounded above by the
isoperimetric profile of the half-space
(Proposition~\ref{prop:firstcomparison}).  Then we prove in
Theorem~\ref{th:half-space}, by using Theorem~\ref{th:convergence} and
a result by P.~B\'erard and D.~Meyer~\cite{bemeyer}, that, in case
of nonexistence of isoperimetric regions for a given volume, the
isoperimetric profile of $M$ is the one of the half-space (see the
notation subsection below for a precise definition).  This is a strong
result that gives us existence of isoperimetric regions under several
conditions on the cone: for instance, when $M$ is a cone over a domain
$C\subset\sph^{n}$, with
$\mathcal{H}^{n}(C)\le\mathcal{H}^{n}(\sph^{n})/2$
(Proposition~\ref{prop:halfvolume}), or when the boundary of the cone
admits a local supporting hyperplane
(Proposition~\ref{prop:existence}).  In particular, for a convex cone,
there is existence for any value of volume.  We also prove that any
isoperimetric region in a cone is bounded in
Proposition~\ref{prop:boundedness}.

In Section~\ref{sec:stability} we show which are the stable regions in
a convex cone $M$.  We define them as regions bounded by hypersurfaces
of constant mean curvature $H$ with a singular set of small Hausdorff
dimension such that the second derivative of perimeter is nonnegative
for any volume preserving variation.  In Theorem~\ref{th:stability} we
prove that the only stable regions in a convex cone are round balls
contained in the closure of the cone, or balls centered at the vertex
intersected with the cone, or half-balls centered and lying over a
flat piece of $\ptl M$.  The method we use to prove
Theorem~\ref{th:stability} was introduced by J.~L.~Barbosa and
M.~do~Carmo in \cite{bdc} to show that round spheres are the only
compact stable constant mean curvature hypersurfaces in $\rr^{n+1}$,
and was adapted later by F.~Morgan and M.~Ritor\'e~\cite{morganritore}
to identify compact stable hypersurfaces of constant mean curvature in
certain cones, allowing the presence of small singular sets.  Once
Theorem~\ref{th:stability} is proved, a simple comparison allows us to
characterize the isoperimetric regions in a convex cone in
Theorem~\ref{th:main}: they are the round balls centered at the vertex
intersected with the cone.  This result was first obtained by
P.~L.~Lions and F.~Pacella \cite{lions-pacella} by using
Brunn-Minkowski theory.  In \cite{morgan-riem}, Remark after
Theorem~10.6, Frank Morgan has pointed that this result could also be
obtained by modifying Gromov's proof of the isoperimetric inequality
in Euclidean space.

We have also included in a final section as an appendix, a direct
proof of Theorem~\ref{th:main} without using stability.  The idea is
as follows: for any isoperimetric region $\Om$ we consider the set of
equidistant hypersurfaces to its boundary; then, for certain function
depending on the relative profile of the equidistants it is satisfied
that the second derivative with respect to the volume is nonnegative. 
From this, we deduce the result after an explicit calculation of the
above derivative.  The possible presence (in high dimensions) of
singularities in the boundary of $\Om$ requires approximation. 
Similar approximation arguments were used by F.~Morgan and
D.~Johnson~\cite{mj} to obtain a differential inequality for the
isoperimetric profile in a compact Riemannian manifold.

\subsection*{Notation and preliminaries}
Let $M^{n+1}$ be an open connected set of a smooth, complete,
$(n+1)$-dimensional Riemannian manifold.  The $(n+1)$-dimensional and
the $k$-dimensional Hausdorff measures of a set $\Om\subseteq M$ will
be denoted by $\vol(\Om)$ and $\mathcal{H}^k(\Om)$ respectively.  For
any measurable set $\Om\subseteq M$ and any open set
$U\subseteq M$, let $\pp(\Om,U)$ be the
\emph{perimeter} of $\Om$ \emph{relative to} $U$, defined as
\[
\pp(\Om,U)=\sup\left\{\int_\Om\text{div}\,Y
\,d\mathcal{H}^{n+1}:|Y|\leq 1\right\},
\]  
where $Y$ is a smooth vector field with compact support contained in
$U$, and $\text{div}\,Y$ is the divergence of $Y$ \cite{chavel}, p.~3. 
\emph{We shall simply denote} $\pp(\Om)=\pp(\Om,M)$.

A set $\Om$ is said to be of \emph{finite perimeter in} $U$ if
$\vol(\Om\cap U)<\infty$ and $\pp(\Om,U)<\infty$.  This condition is
equivalent to requiring that the characteristic function of $\Om$ is
of \emph{bounded variation} in $U$.  In case $U=M$ we simply say that
$\Om$ is of finite perimeter.  For further background about perimeter
and sets of finite perimeter we refer to the reader to \cite{giusti}
and \cite{bvfunctions}.

As the perimeter of a set is not changed by adding or removing sets of
$\mathcal{H}^{n+1}$-measure $0$, we shall always assume that
$0<\vol(\Om\cap B_{r}(p))<\vol(B_{r}(p))$, for every ball $B_{r}(p)$
centered at $p\in\ptl\Om$.

The \emph{isoperimetric profile} of $M$ is the function $I_M:(0,\vol
M)\to\rr^{+}\cup\{0\}$ given by
\[
I_M(V)=\inf\{\pp(\Om):\,\Om\subseteq M,\,\, \vol(\Om)=V\}.
\]
An \emph{isoperimetric region} in $M$ for volume $V\in (0,\vol M)$ is
a set $\Om\subseteq M$ satisfying $\vol(\Om)=V$ and $\pp(\Om)=I_M(V)$. 
A {\it minimizing sequence\/} of sets of volume $V$ is a sequence of
sets of finite perimeter $\{\Om_{k}\}_{k\in\nn}$ such that
$\vol(\Om_k)=V$ for all $k\in\nn$ and
$\lim_{k\to\infty}\pp(\Om_k)=I_M(V)$.  We recall that a sequence
$\{\Om_k\}_{k\in\nn}$ converges in the finite perimeter sense to a set
$\Om$ if $\{f_k\}\to f$ in $L^1_{\text{loc}}(M)$ and
$\lim_{k\to\infty}\pp(\Om_k)=\pp(\Om)$.  Here, $f_k$ and $f$ denote
the characteristic function of $\Om_k$ and $\Om$ respectively.

Finally, we remark the following consequence of the coarea formula
\cite{burago}, Theorem~13.4.2; for any measurable set $\Om\subset M$,
we can compute $\vol(\Om)$ as
\begin{equation*}
\vol(\Om)=\int_0^{+\infty}\mathcal{H}^n(\Om\cap S_t)\,dt,
\end{equation*}
where $S_t$ is the metric sphere in $M$ of radius $t$ centered at a
point $p_0\in M$.

\section{Minimizing sequences and regularity of isoperimetric regions
in a Riemannian manifold}
\label{sec:existence}
\setcounter{equation}{0}

In this section we study the behavior of a minimizing sequence for
fixed volume.  We shall work in open sets of complete Riemannian
manifolds.  The perimeter will be the one relative to the open set. 
The technique can be extended, with no extra effort, to other
situations.

\begin{theorem}
\label{th:convergence}
Let $M^{n+1}$ be a connected non bounded open set of a
complete~Riemannian manifold.  For any minimizing sequence
$\{\Om_k\}_{k\in\nn}$ of sets of volume $V$, there exist a finite
perimeter set $\Om\subset M$ and sequences of sets of finite perimeter
$\{\Om_k^c\}_{k\in\nn}$, $\{\Om_k^d\}_{k\in\nn}$, such that
\begin{enum} 
\item $\vol(\Om)\leq V,\quad\pp(\Om)\leq I_M(V)$.  
\item $\vol(\Om_k^c)+\vol(\Om_k^d)=V,\quad
\lim_{k\to\infty}[\pp(\Om_k^c)+\pp(\Om_k^d)]=I_M(V)$.
\item The sequence $\{\Om_k^d\}_{k\in\nn}$ diverges.
\item Passing to a subsequence, we have
$\lim_{k\to\infty}\vol(\Om_k^c)=\vol(\Om)$ and
$\lim_{k\to\infty}\pp(\Om_k^c)=\pp(\Om)$.  In fact,
$\{\Om_k^c\}_{k\in\nn}$ converges to $\Om$ in the finite perimeter
sense.
\item $\Om$ is an isoperimetric region \emph{(}it could be
empty\emph{)} for the volume it encloses.
\end{enum}
\end{theorem}

\begin{remark}
For $n=1$ it was shown in \cite{ritore2}, Lemma~2.2, that the sets
$\Om_{k}^{c}$ and $\Om_{k}^{d}$ can be taken as union of connected
components of $\Om_{k}$.
\end{remark}

\begin{proof}
Fix a point $p_0\in M$ and denote by $B_r$ (resp.~$S_r$) the
intersection of $M$ with the ball (resp.~sphere) in the ambient
manifold centered at $p_{0}$ of radius $r>0$.  Clearly $B_m\subset
B_{m+1}$ for $m\in\nn$ and $M=\bigcup_{m\in\nn}\,B_m$.

Let $f_k$ be the characteristic function of $\Om_k$.  Applying the
Compactness Theorem~\cite{giusti}, Theorem 1.19, on each $B_m$ and a
diagonal argument, we get the existence of a function $f\in
L^1_{\text{loc}}(M)$ and of a subsequence, that we denote again by
$\{\Om_k\}$, such that $\{f_k\}\to f$ in $L^1_{\text{loc}}(M)$.  We
can even assume that $\{f_k\}$ pointwise converges
$\mathcal{H}^{n+1}$-almost everywhere to $f$, and so $f$ is the
characteristic function of a set $\Om\subset M$.  By Fatou's~Lemma and
the lower semicontinuity of perimeter \cite{giusti}, Theorem 1.9, we
have $\vol(\Om)\leq V$ and $\pp(\Om)\leq I_M(V)$.  If $\vol(\Om)=V$,
then $\Om$ is an isoperimetric region of volume $V$ and the conclusion
follows by taking $\Om_k^c=\Om_k$ and $\Om^d_k=\emptyset$ for all
$k\in\nn$.  Henceforth we assume that $\vol(\Om)<V$.

To define $\Om_k^c$ and $\Om^d_k$ we need to choose an appropriate
truncation of $\Om_k$.  By the coarea formula, on any real interval
$J$ we have $\int_{J} \mathcal{H}^{n}(\Om_{k}\cap
S_{t})\,dt\le\text{vol}(\Om_{k})=V$.  If $\text{length}(J)\ge\ell$
then there exists $t\in J$ such that
\[
\mathcal{H}^{n}(\Om_{k}\cap S_{t})\le\frac{V}{\ell}.
\]
From this observation, passing again to a subsequence of $\Om_{k}$ if
necessary, it is easy to find a sequence $r(k)$ of positive real 
numbers such that
\begin{align}
\label{eq:length}
r(k+1)-r(k)&\ge k, \\
\label{eq:hausdorff}
\mathcal{H}^{n}(\Om_{k}\cap S_{r(k+1)})&\le\frac{V}{k}, \\
\label{eq:integral}
\int_{B_{r(k+1})} |f_{k}-f|\,d\mathcal{H}^{n+1}&\le\frac{1}{k}.
\end{align}

For each $k\in\mathbb{N}$, define $\Om^c_k=\Om_k\cap B_{r(k+1)}$ and
$\Om^d_k=\Om_k-B_{r(k+1)}$.  Clearly $\Om_k=\Om^c_k\cup\Om^d_k$ and so
$\vol(\Om_k^c)+\vol(\Om^d_k)=V$.  Moreover, we can choose the sequence
$\{r(k)\}_{k\in\nn}$ so that \cite{bvfunctions}, Corollary 5.5.3,
\begin{align*}
\pp(\Om_k^c)&\leq\pp(\Om_k, B_{r(k+1)})+\mathcal{H}^n(\Om_k\cap S_{r(k+1)}),  
\\
\pp(\Om^d_k)&\leq\pp(\Om_k,M-\overline{B}_{r(k+1)})+
\mathcal{H}^n(\Om_k\cap S_{r(k+1)}).
\end{align*}

\noindent
As $\text{vol}(\Om_k)=V$, we get from \eqref{eq:hausdorff}
\[
I_M(V)\leq\pp(\Om_k)\leq\pp(\Om_k^c)+
\pp(\Om^d_k)\leq \pp(\Om_k)+
2\,\mathcal{H}^n(\Om_k\cap S_{r(k+1)})\leq\pp(\Om_k)+\frac{2V}{k},
\]
which implies that $\lim_{k\to\infty}\,[\pp(\Om_k^c)+\pp(\Om^d_k)]=I_M(V)$ 
and proves (ii). Assertion (iii) follows from the inclusion 
$\Om^d_k\subset M-B_{r(k+1)}$.

To prove (iv) observe
\begin{align*}
|\vol(\Om_k^c)-\vol(\Om)|&\leq|\vol(\Om_k\cap B_{r(k+1)})-
\vol(\Om\cap B_{r(k+1)})|+\vol(\Om- B_{r(k+1)})  
\\
&\leq\int_{B_{r(k+1)}}|f_k-f|\,\,d\mathcal{H}^{n+1}+\vol(\Om-B_{r(k+1)}),  
\end{align*}
and so $\lim_{k\to\infty}\vol(\Om_k^c)=\vol(\Om)$ by
\eqref{eq:integral}.  On the other hand, by the definition of the sets
$\Om_k^c$ we have
\[
\int_{B_{r(k+1)}}|f^c_t-f|\,\,d\mathcal{H}^{n+1}=
\int_{B_{r(k+1)}}|f_t-f|\,\,d\mathcal{H}^{n+1},\quad\ t\geq k,
\] 
where $f^c_t$ denotes the characteristic function of $\Om_t^c$.  The
above equality gives us $\{f^c_k\}_k\to f$ in $L^1_{\text{loc}}(M)$. 
Hence $\pp(\Om)\leq\liminf_{k\to\infty}\pp(\Om_k^c)$ by the lower
semicontinuity of perimeter.

To finish the proof we only have to see that $\Om$ is
perimeter minimizing among sets of volume $V_{0}=\vol(\Om)$, and
that $\liminf_{k\to\infty}\pp(\Om_{k}^c)=\pp (\Om)$.  If any of these 
statements is not true, then we will construct an ``improved 
minimizing sequence'' to get a contradiction.

Suppose first that $\Om$ is not an isoperimetric region.  In this
case, we can find a set $E\subset M$ of volume $V_{0}$ with smooth
boundary such that $\pp(E)<\pp(\Om)$.  From \eqref{eq:length} there 
exists $t(k)\in (r(k),r(k+1))$ such that
\[
\mathcal{H}^n(E\cap S_{t(k)})\le\frac{V_0}{k}.
\]
Let $E_{t(k)}=E\cap B_{t(k)}$.  Obviously
$\lim_{k\to\infty}\vol(E_{t(k)})=V_0$.

Now, take a point $q_0\in\ptl E\cap M$ and a small neighborhood $W$ of
$q_0$ in $\ptl E$.  Let $u\not\equiv 0$ be a nonnegative function in
$C^\infty_0(W)$.  For all $k$ large enough $W\subset E_{t(k)}$ and $u$
induces a variation of $\ptl E_{t(k)}$ that can be used to modify
slightly the volume of $E_{t(k)}$.  Since $\text{vol}(E_{t(k)})\to
V_{0}$ and $\text{vol}(\Om_{k}^{d})\to V-V_{0}$, we use the variation
associated to $u$ to get, for large $k$, a set $E_{t(k)}'$ so that
\[
\text{vol}(E_{t(k)}')+\text{vol}(\Om_{k}^{d})=V.
\]
By the first variation for perimeter and volume, we have
\begin{equation}
\label{eq:perimetertk}
\pp(E'_{t(k)})\leq\pp(E_{t(k)})+\alpha\,|v_{k}|,
\end{equation}
where $v_{k}=\text{vol}(E_{t(k)}')-\text{vol}(E_{t(k)})$ converges to 
$0$ and $\alpha$ is a positive constant independent of~$k$. 

Let $T_k=E'_{t(k)}\cup\Om_k^d$.  The volume of $T_k$ equals $V$ and, by
\eqref{eq:perimetertk},
\begin{align*}
\pp(T_k)&\leq\pp(E,B_{t(k)})+
\mathcal{H}^n(E\cap S_{t(k)})+\alpha\, |v_k|+
\pp(\Om_{k}^d)
\\
&\leq\pp(E)+\frac{V_0}{k}+\alpha\, |v_k|+\pp(\Om_{k}^d),  
\end{align*}
from which we obtain 
\[
\liminf_{k\to\infty}\pp(T_k)\leq\pp(E)+
\liminf_{k\to\infty}\pp(\Om_k^d)<\pp(\Om)+
\liminf_{k\to\infty}\pp(\Om_k^d)\leq I_M(V),
\]
since $\pp(\Om)\leq\liminf_{k\to\infty}\pp(\Om_k^c)$.  By this
contradiction (v) is proved.  The proof of (iv) is similar; in fact we
only have to work with $\Om$ instead of $E$.  We must remark that this
can be done since the set of regular points in $\ptl\Om\cap M$ is
open, see Proposition~\ref{prop:regularity} below.
\end{proof}

We finish this section by recalling smoothness properties of the
boundary of an isoperimetric region in an open subset $M^{n+1}$ of an
$(n+1)$-dimensional Riemannian manifold.  Regularity results for
isoperimetric hypersurfaces in the Euclidean setting were obtained by
E.~Gonzalez, U.~Massari and I.~Tamanini~\cite{gmt}, who treated
interior regularity, and by M.~Gr\"{u}ter~\cite{gruter}, who studied
regularity near boundary points.  As pointed out by F.
Morgan~\cite{morgan} their results can also be stated in the setting
of Riemannian manifolds, by using the paper of
F.~Almgren~\cite{almgren}.  A complete proof of interior
regularity for isoperimetric boundaries in Riemannian manifolds can be
found in \cite{morgan2}; weaker regularity properties when $\ptl M$ is
not $C^\infty$, but only $C^{1,1}$, are also established.  We collect
the above results in next Proposition.

\begin{proposition}[{{\bf Regularity}}]
\label{prop:regularity}
Let $\Om$ be a perimeter-minimizing set under a volume cons\-traint in
a connected open set $M^{n+1}$ with smooth boundary $\ptl M$ of an
$(n+1)$-dimensional Riemannian manifold.  Then, the boundary
$\Lambda=\overline{\ptl\Om\cap M}$ can be written as a disjoint union
$\Sg\cup\Sg_0$, where $\Sg$ is the regular part of $\Lambda$ and
$\Sg_0=\Lambda-\Sg$ is the set of singularities.  Moreover, we have
\begin{enum}
\item $\Sg\cap M$ is a smooth, embedded hypersurface with constant
mean curvature with respect to the inner normal.  \item If
$x\in\Sg\cap\ptl M$, then $\Sg$ is a smooth, embedded hypersurface
with boundary contained in $\ptl M$ in a neighborhood of $x$; in this
neighborhood $\Sg$ has constant mean curvature and meets $\ptl M$
orthogonally.
\item $\Sg_0$ is a closed set of Hausdorff dimension less than or
equal to $n-7$.
\end{enum}
\end{proposition}

\begin{remark}
In the preceding Proposition the regular set $\Sg$ is defined as
follows: for $x\in\Sg$ there is a neighborhood $W$ of $x$ in $\Sg$
such that $W$ is a smooth, embedded hypersurface without boundary or
with boundary contained in $\ptl M$.  Note that one of the conclusions
of the above Proposition is the absence of interior points in $\Sg$
meeting $\ptl M$ tangentially, see \cite{gruter-nuevo}.
\end{remark}

\begin{example}
This example illustrates that isoperimetric regions in a smooth
manifold with boundary need not meet the boundary of the manifold.
 
Let $M$ be the compact surface obtained from attaching the hemisphere
of $\mathbb{S}^2$ centered at the north pole to the compact cylinder
$\mathbb{S}^1\times [-h,0]$ ($h>\pi$) through the circle
$\mathbb{S}^1\times\{0\}$.  Suppose that $\Om$ is an isoperimetric
region in $M$ of area $A>2\pi$ such that
$\Lambda=\overline{\ptl\Om\cap M}$ meets $\ptl M$ orthogonally.  As
the cylinder is locally isometric to the half-plane at any point of
$\ptl M$, we deduce that $\Lambda$ contains either a semicircle
centered at $\ptl M$ or two vertical segments $\{\theta_i\}\times
[-h,0]$, $i=1,2$.  In any of these cases we have
$\pp(\Om)^2>4\pi^2=\pp(D)^2$, where $D$ is the geodesic disk in $M$
centered at the north pole and enclosing area $A$.  This contradiction
shows that any isoperimetric region $\Om$ of area $A>2\pi$ satisfies
$\Lambda\cap\ptl M=\emptyset$.  If $\overline{\Om}\cap\ptl M$ is
nonvoid, then $\Om'=M-\Om$ is an isoperimetric region that does not
meet $\ptl M$.
\end{example}

\section{Existence of isoperimetric regions inside a Euclidean cone}
\label{sec:cone}
\setcounter{equation}{0}

In this section we shall consider a cone $M=0\cone C$, where $C$ is a
connected open set (domain) in $\sph^n\subset\rr^{n+1}$ with smooth
boundary.  The closure of the cone is $\overline{M}=0\cone
\overline{C}$.  When $C$ is an open half-sphere the cone coincides with
an open half-space.  It is well known that the isoperimetric regions in
a half-space are half-balls centered at the boundary of the half-space
\cite{zumbrun1}, and that the isoperimetric profile of the half-space is
given by
\[
I_n(V)=\left (\frac{\mathbf{c}_n}{2}\right )^{1/(n+1)}
(n+1)^{n/(n+1)}\, V^{n/(n+1)},\qquad V>0,
\]
where $\mathbf{c}_n=\mathcal{H}^n(\mathbb{S}^n)$.

We begin with two easy consequences of the invariance of a cone by
dilations.

\begin{proposition}
\label{prop:profile}
The isoperimetric profile $I_M$ of a cone $M$ satisfies
\begin{equation*}
I_M(V)=c\,V^{n/(n+1)}\quad\text{ for } V>0,
\end{equation*}
where $c$ is a constant equal to $I_M(V_{0})/V_{0}^{n/(n+1)}$ for any
$V_{0}>0$.  Hence $I_M^{(n+1)/n}(V)$ is a linear function of $V$.
\end{proposition}

\begin{proof}
For any $\lambda>0$ let $h_{\lambda}(p)=\lambda p$.  Fix $V_{0}>0$ and
consider a finite perimeter set $\Om\subset M$ with $\vol(\Om)=V_0$. 
For each $\lambda>0$, the set $\Om_\lambda=h_\lambda(\Om)$ is
contained in $M$ and satisfies $\vol(\Om_\lambda)=\lambda^{n+1}V_0$
and $\pp(\Om_\lambda)=\lambda^n\,\pp(\Om)$.  By the definition of
$I_M$ we have $I_M(\lambda^{n+1}V_0)\leq\lambda^n\,\pp(\Om)$, and so
$I_M(\lambda^{n+1}V_0)\leq\lambda^n\, I_M(V_0)$.  With a similar
argument we obtain the opposite inequality and
\begin{equation*}
\lambda^n I_M(V_{0})=I_M(\lambda^{n+1}V_{0}).
\end{equation*}
From this equality the result follows by taking
$\lambda=(V/V_0)^{1/(n+1)}$.
\end{proof}

\begin{proposition}
\label{prop:exornot}
Let $M$ be a cone over a smooth domain in a sphere
$\sph^{n}\subset\rr^{n+1}$.  Then either there exist isoperimetric
regions for any given value of volume, or there are no isoperimetric
regions at all.
\end{proposition}

\begin{proof}
Let $h_\lambda(p)=\lambda p$.  It is clear that $\Om\subset M$ is an
isoperimetric region of volume $V_{0}$ if and only if
$h_{\lambda}(\Om)\subset M$ is an isoperimetric region of volume
$\lambda^{n+1}V_{0}$.
\end{proof}

The following result gives an upper bound on the isoperimetric profile
of a cone $M$ and it is a consequence of the fact that $\ptl M$ is
asymptotically flat.

\begin{proposition}
\label{prop:firstcomparison}
The isoperimetric profile $I_M$ of a cone $M\subset\rr^{n+1}$ is
bounded from above by the isoperimetric profile of the half-space.
\end{proposition}

\begin{proof}
Consider a diverging sequence of metric balls $D_k'$ centered at
points of $\ptl M$, so that the volume of $D_k=D_k'\cap M$ equals $V$. 
By the definition of isoperimetric profile we have
$I_M(V)\leq\pp(D_k)$ for all $k\in\nn$.  On the other hand
$\lim_{k\to\infty}\pp(D_k)=I_n(V)$ since $\ptl M$ approaches a
hyperplane at infinity.  Thus,
\[
I_M(V)\leq\lim_{k\to\infty}\pp(D_k)=I_n(V).
\]
\end{proof}

The main result of this section is the following Theorem

\begin{theorem}
\label{th:half-space}
Let $M$ be a cone in $\rr^{n+1}$.  If there are not isoperimetric
regions in $M$, then the isoperimetric profile $I_M$ of $M$ is the one
of the half-space.
\end{theorem}

\begin{proof}
We know from Proposition~\ref{prop:firstcomparison} that inequality
$I_M(V)\leq I_n(V)$ holds for all $V>0$.  To see the reverse
inequality suppose that there are not isoperimetric regions in $M$. 
From Theorem~\ref{th:convergence} this means that the volume of the
sets of the convergent part of any minimizing sequence of sets of
volume $V$ goes to zero.  Hence there exists a diverging sequence
$\{\Om_k\}_{k\in\nn}$ of sets of finite perimeter in $M$ such that
\begin{enum}
\item $\vol(\Om_k)\leq V$ for all $k\in\nn$.
\item $\lim_{k\to\infty}\vol(\Om_k)=V$.
\item $\lim_{k\to\infty}\pp(\Om_k)=I_M(V)$. 
\end{enum}

Fix an arbitrary $\eps>0$.  By B\'erard-Meyer isoperimetric inequality
\cite{bemeyer}, App.~C, there exists a constant $V_\eps>0$ such that,
for any set of finite perimeter $\Om\subset B_{4}-\overline{B}_{1}$
with $\vol(\Om)\leq V_\eps$, we get
\begin{equation}
\label{eq:berardmeyer}
\pp(\Om, B_4-\overline{B}_1)^{(n+1)/n}\geq c_n\,
(1-\eps)^{(n+1)/n}\,\vol(\Om),
\end{equation}
where $B_r$ is the intersection of the open ball of radius $r$
centered at the origin with the cone, and $c_n$ is the constant that
appears in the expression of the isoperimetric profile of the
half-space $I_n(V)^{(n+1)/n}=c_{n}\,V$.  By scaling we get that, for
any finite perimeter set $\Om\subset
B_{4\lambda}-\overline{B}_{\lambda}$ such that
$\vol(\Om)\leq\lambda^{n+1}V_\eps$, we have
\begin{equation}
\label{eq:berardmeyer2}
\pp(\Om, B_{4\lambda}-\overline{B}_{\lambda})^{(n+1)/n}\geq c_n\,
(1-\eps)^{(n+1)/n}\,\vol(\Om).
\end{equation}

Let $S_{r}$ be the intersection of the sphere of radius $r>0$ centered
at the vertex of the cone with the cone.

Since the sequence $\{\Om_k\}_k$ diverges we may assume, passing to a
subsequence if necessary, that $\overline{\Om}_{m}\subset
M-\overline{B}_m$.

Fix $m\in\nn$ so that $m^{n+1}\ge V/V_{\eps}$.  For such $m$ we use
the coarea formula to define a sequence $\{t(q)\}_{q\in\nn}$ of
positive real numbers by setting $t(0)=m$ and $t(q)\in
(2^{q-1}m,2^{q}m)$, $q\ge 1$, so that
\begin{equation}
\label{eq:bm2}
\mathcal{H}^n(\Om_{m}\cap S_{t(q)})\leq
\frac{\vol(\Om_{m})}{2^{q-1}\,m}\leq\frac{V}{2^{q-1}\,m}, \qquad q\ge 
1.
\end{equation}
Obviously $t(1)/t(0)\le 2$ and $t(q)/t(q-1)< 4$ for $q\ge 2$.

Define $\Om^q_{m}=\Om_{m}\cap (B_{t(q)}-\overline{B}_{t(q-1)})$. The 
set $\Om_{m}^{q}$ is contained in $B_{t(q)}-\overline{B}_{t(q-1)}$. We 
have $t(q)< 4\,t(q-1)$, and also
\[
\text{vol}(\Om_{m}^{q})\le\text{vol}(\Om_{m})\le V\le m^{n+1}V_{\eps}\le
t(q-1)^{n+1}V_{\eps}.
\]
From inequality~\eqref{eq:berardmeyer2} we get
\[
\pp(\Om^q_{m}, B_{4\,t(q-1)}-\overline{B}_{t(q-1)})^{(n+1)/n}\geq
c_n\,(1-\eps)^{(n+1)/n}\,\vol(\Om^q_{m}),\qquad q\in\nn,
\]
and so, from inequality $\sum_{i} a_{i}^{(n+1)/n}\le (\sum_{i} 
a_{i})^{(n+1)/n}$,
\begin{equation}
\label{eq:sum}
\left[\sum_{q=1}^\infty\pp(\Om_{m}^q,
B_{4\,t(q-1)}-\overline{B}_{t(q-1)})\right]^{(n+1)/n}\geq
c_n\,(1-\eps)^{(n+1)/n}\,\vol(\Om_{m}).
\end{equation}
Taking into account that
\[
\pp(\Om_{m}^q, B_{4\,t(q-1)}-\overline{B}_{t(q-1)})\le
\pp(\Om_{m}^{q}, B_{t(q)}-\overline{B}_{t(q-1)})+
\mathcal{H}^{n}(\Om_{m}\cap S_{t(q)}), \qquad m\in\nn,
\]
and that
\[
\sum_{q=1}^\infty\pp(\Om_{m}^q, B_{t(q)}-\overline{B}_{t(q-1)})
=
\sum_{q=1}^\infty\pp(\Om_{m}, B_{t(q)}-\overline{B}_{t(q-1)})\le
\pp(\Om_{m}),
\]
we have, from \eqref{eq:bm2},
\[
\pp(\Om_{m})+\sum_{q=1}^{\infty} \frac{V}{2^{q-1}m}\ge
\sum_{q=1}^{\infty} \pp(\Om_{m}^{q}, 
B_{4\,t(q-1)}-\overline{B}_{t(q-1)}).
\]

From this inequality and \eqref{eq:sum}, letting $m\to\infty$ we
conclude
\[
I_M(V)^{(n+1)/n}\geq c_n\,(1-\eps)^{(n+1)/n}\,V=(1-\eps)^{(n+1)/n}\,
I_n(V)^{(n+1)/n}.
\]
As $\eps$ is  arbitrary the proof follows.
\end{proof}

Now, we establish some criteria to decide whether there exist or not
isoperimetric regions in a cone.

\begin{proposition}
\label{prop:halfvolume}
Let $M\subset\rr^{n+1}$ be a cone over a domain $C\subset\sph^n$ such
that $\mathcal{H}^n(C)\leq\mathbf{c}_{n}/2$, where $\mathbf{c}_n$ is
the Riemannian volume of\/ $\sph^n$.  Then isoperimetric regions exist
in $M$ for any value of volume.

In particular, isoperimetric regions of any volume exist in a convex
cone.
\end{proposition}

\begin{proof}
Simply observe that the perimeter and the volume of a ball $B_r$ of
radius $r>0$ centered at the origin intersected with $M$ are given by
\[
\pp(B_{r})=r^n\,\mathcal{H}^n(C),
\qquad\text{vol}(B_r)=\frac{r^{n+1}}{n+1}\,\mathcal{H}^n(C).
\]
The relative isoperimetric profile of these balls equals
\begin{equation*}
\pp(V)=\mathcal{H}^n(C)^{1/(n+1)}\,(n+1)^{n/(n+1)}\,V^{n/(n+1)},
\end{equation*}
which is less than or equal to the isoperimetric profile of the
half-space.  Hence the first part of the Proposition follows by
Theorem~\ref{th:half-space}.

If $M$ is convex then $\mathcal{H}^n(C)\le\mathbf{c}_{n}/2$ and we 
apply the first part of the Proposition.
\end{proof}

Another criterion for existence of isoperimetric regions in a cone is
given in next Proposition.

\begin{proposition}
\label{prop:existence}
Consider a cone $M$ over a smooth domain $C\subset\sph^n$, and assume
that $\ptl M-{\{0}\}$ has a point $x$ admitting a local supporting
hyperplane.  Then isoperimetric regions exist in $M$ for any value of
volume.
\end{proposition}

\begin{proof}
Let $\pp(r)$ and $V(r)$ be the perimeter and the volume of the ball
$B_r$ of radius $r>0$ centered at $x$ intersected with the cone, and
$\wt{V}(r)$ the volume of the cone subtended by $\ptl B_{r}\cap M$ and
vertex at $x$.  We have the relation
\[
\pp(r)=(n+1)\,\frac{\wt{V}(r)}{r}.
\]
Since $M$ is locally convex near $x$, $V(r)\ge\wt{V}(r)$ for $r$
small, so that
\[
\pp(r)=(n+1)\,\frac{\wt{V}(r)}{r}\le (n+1)\,\frac{V(r)}{r}.
\]
On the other hand, if $\pp_{e}(r)$ and $V_{e}(r)$ are the area and 
volume of the half-ball of radius $r>0$, we have
\[
\frac{\pp_{e}(r)}{V_{e}(r)}=\frac{n+1}{r},
\]
and so
\[
\frac{\pp(r)}{V(r)}\le\frac{\pp_{e}(r)}{V_{e}(r)}.
\]
Since $V(r)\le V_{e}(r)$ by the local convexity of $\ptl M$ around $x$
we get
\[
\frac{\pp(r)}{V(r)^{n/(n+1)}}=\frac{\pp(r)}{V(r)}\,V(r)^{1/(n+1)}\le
\frac{\pp_{e}(r)}{V_{e}(r)}\,V_{e}(r)^{1/(n+1)}=
\frac{\pp_{e}(r)}{V_{e}(r)^{n/(n+1)}}=d_{n},
\]
where $d_{n}$ is the constant that appears in the expression of the 
isoperimetric profile of the half-space $I_{n}(V)=d_{n}V^{n/(n+1)}$.

Hence, for small $r$, $\mathcal{P}({B_r})\le I_{n}(V(r))$.  By
Theorem~\ref{th:half-space} isoperimetric regions exist in the cone for
any value of volume.
\end{proof}

The last result of this section shows the boundedness of an
isoperimetric region in a cone.  The proof is modelled on the
Euclidean one \cite{morgan}, Lemma 13.6.  We include it here for the
sake of completeness.

\begin{proposition}
\label{prop:boundedness}
Any isoperimetric region $\Om$ in a cone $M$ is bounded.
\end{proposition}

\begin{proof}
For any $r\geq 0$ denote $V(r)=\vol(\Om-B_r)$ and $\pp(r)=\pp(\Om,
M-\overline{B}_r)$, where $B_r$ is the open ball of radius $r$
centered at the origin.  The function $V(r)$ is decreasing and
satisfies $\lim_{r\to 0}V(r)=0$.  By the coarea formula we have
$V'(r)=-\mathcal{H}^n (\Om\cap\partial B_r)$ for almost all $r\geq 0$.

By applying the isoperimetric inequality of the cone in
Proposition~\ref{prop:profile} and
\cite{bvfunctions}, Corollary~5.5.3, we obtain
\begin{equation}
\label{eq:ones}
|V'(r)|+\pp(r)\ge\pp(\Om-B_r)\geq c\,V(r)^{n/(n+1)},
\end{equation} 
for almost all $r\geq 0$.

Now, take a large $r\geq 0$ so that $V(r)$ is close enough to zero. 
By (i) in Proposition~\ref{prop:regularity} we can select an open
subset $W$ of regular points of $\ptl\Om\cap M$ contained in
$\Om\cap B_r$, and a negative smooth function $u$ with compact support
contained in $W$.  This function corresponds to a normal variation of
$W$ and, by the first variation formulae for perimeter and volume,
$d\pp/dV(\vol(\Om\cap B_r))=nH$ ($H$ is the constant mean curvature of
$W$).  Hence, if we produce a small increment of volume $V(r)$, then
the increment of perimeter is less than or equal to $(1+nH)V(r)$.  By
using that $\Om$ is an isoperimetric region we deduce
\[
\pp(\Om,B_r)+\mathcal{H}^n(\Om\cap\partial B_r)+(1+nH) \,V(r)\geq
\pp(\Om)\ge\pp(\Om, B_r)+\pp(r),
\] 
which implies that
\begin{equation}
|V'(r)|+(1+n H)\,V(r)\geq\pp(r),\quad\text{ for almost all  }r\geq 0.
\label{eq:two}
\end{equation}

Adding inequalities \eqref{eq:ones} and \eqref{eq:two} and using that
$\lim_{r\to+\infty}V(r)=0$, we obtain
\[
2\,|V'(r)|\geq \frac{c}{2}\,V(r)^{n/(n+1)}, \quad\text{ for almost all
}r\gg 0.
\]

Now, if we admit that $\Om$ is unbounded, then $V(r)>0$ and
\[
(n+1)\,(V^{1/(n+1)})'(r)=-|V'(r)|\,V(r)^{-n/(n+1)}\leq-\,\frac{c}{4},
\]
which is clearly a contradiction since $V(r)$ is positive for all $r>0$.
\end{proof}

\section{Characterization of stable regions with small singular set in
convex cones}
\label{sec:stability}
\setcounter{equation}{0}

We begin this section with a general notion of stability for sets of
finite perimeter.  Let us consider a cone $M\subset\rr^{n+1}$ and a
bounded set of finite perimeter $\Om\subset M$.  We take a vector
field $X$ in $\rr^{n+1}$ with compact support such that $X(p)\in
T_{p}(\ptl M)$ for all $p\in \ptl M-\{0\}$ and we consider the
associated one-parameter group of diffeomorphisms
$\{\varphi_{t}\}_{t}$ of $\overline{M}$.

We shall say that $\Om_{t}=\varphi_{t}(\Om)$ is a
\emph{volume-preserving deformation} of $\Om$ if
$\vol(\Om_t)=\vol(\Om)$ for $t$ small enough.  Let
$\mathcal{P}(t)=\mathcal{P}(\Om_{t})$.  It is said that $\Om$ is
\emph{stationary} if $\pp'(0)=0$ for all volume-preserving
deformations.  It is said that $\Om$ is \emph{stable} if it is
stationary and $\pp''(0)\geq 0$ for all volume-preserving
deformations.  Obviously any isoperimetric region in $M$ is a stable
region.  Theorems~\ref{th:stability} and \ref{th:main} show that, in
general, there are stable regions inside a cone which are not
isoperimetric.

From now on we assume that $\Lambda=\overline{\ptl\Om\cap M}$ can be
decomposed as $\Sg\cup\Sg_0$, where $\Sg_0$ is a closed singular set
consisting of isolated points or such that
$\mathcal{H}^{n-2}(\Sg_0)=0$, and $\Sg$ is a smooth, embedded
hypersurface such that $\ptl\Sg-\Sg_{0}\subset\ptl M$.  We remark that
we are not excluding the possibility that $\Sg$ and $\ptl M$ are
tangent at some interior point to $\Sg$. We define 
$\ptl^{*}\Sg=\ptl\Sg-\Sg_{0}$. If $\ptl^{*}\Sg=\emptyset$, we adopt 
the convention that the integrals over this set are all equal to $0$.

The above setting allows us to give an analytical expression for
stationarity and stability from the variational formulae for perimeter
and volume.  If $X$ is a vector field over $M$ with compact support on
the regular boundary $\Sg$ and such that $X(p)\in T_p(\ptl M)$
whenever $p\in\ptl M-\{0\}$, then $X$ induces a variation of $\Om$ for
which the derivatives of perimeter and volume are given by
\[
\pp'(0)=-\int_{\Sg}nH\,u\,d\mathcal{H}^n-
\int_{\ptl^{*}\Sg} \escpr{X,\nu}\,d\mathcal{H}^{n-1},\qquad
V'(0)=-\int_\Sg u\,d\mathcal{H}^n,
\]
where $H$ is the mean curvature with respect to a unit normal vector
field $N$, $u=\escpr{X,N}$ is the normal component of the variation,
and $\nu$ is the inward normal vector to $\ptl^{*}\Sg$ in $\Sg$. 
Hence, if $\Om$ is stationary, then $H$ is constant and $\Sg$ meets
$\ptl M$ orthogonally in the points of $\ptl^{*}\Sg$, \cite{zumbrun1}.

The derivative $\pp''(0)$ of perimeter when $\Om$ is stationary and
the variation preserves volume is given by \cite{rosvergasta}
\begin{equation}
\label{eq:indexform}
Q(u,u)=-\int_{\Sg} u\,(\Delta u+|\sg|^2u)\,d\mathcal{H}^n
-\int_{\ptl^{*}\Sg}u\,\left\{\frac{\ptl u}{\ptl\nu}
+\text{II}(N,N)\,u\right\}d\mathcal{H}^{n-1},
\end{equation}
where $\Delta$ is the Laplacian on $\Sg$, $|\sigma|^2$ is the squared
sum of the principal curvatures associated to $N$, and $\text{II}$ is
the second fundamental form of $\ptl M-\{0\}$ with respect to the
inner normal.  The above equation defines a quadratic form on
$C^\infty_0(\Sg)$ called the \emph{index form}.

As in \cite{bdc}, a function $u\in C^\infty_0(\Sg)$ satisfying
$\int_\Sg u\,d\mathcal{H}^n=0$ induces a volume-preserving deformation
of $\Om$ with associated vector field $X$ so that the normal component
$\escpr{X,N}$ equals $u$.  Thus, if $\Om$ is stable, then we have
$Q(u,u)\geq 0$ for any function $u$ with compact support and mean zero
on $\Sg$.  We conclude that if $\Om$ is a stable region in $M$ then
$\Sg$ is a stable constant mean curvature hypersurface as defined in
\cite{bdc} and \cite{rosvergasta}.  Moreover, for a function $u\in
C^{\infty}_{0}(\Sg)$, integration by parts gives another expression
for the index form, namely
\begin{equation}
\label{eq:indexform2}
I(u,u)=\int_\Sg (|\nabla u|^2-|\sg|^2u^2)\,d\mathcal{H}^n-
\int_{\ptl^{*}\Sg} \text{II}(N,N)\,u^2\,d\mathcal{H}^{n-1},
\end{equation}
where $\nabla u$ is the gradient of $u$ on $\Sg$.  The right side of
\eqref{eq:indexform2} has meaning for functions in the Sobolev space
$H^{1}(\Sg)$ with compact support in $\Sg$.

Inequality $Q(u,u)\geq 0$ gives interesting geometrical and
topological information when a suitable function $u$ is inserted. 
Following \cite{bdc} we shall prove Theorem~\ref{th:stability} by
considering the test function $u=1+H\escpr{X,N}$, where $X$ is the
position vector field in $\rr^{n+1}$ given by $X(p)=p$.  It was proved
by H.~Wente~\cite{wente} that the function $u$ appears when one
consider first a contraction of a smooth hypersurface by parallels and
then apply a dilation to restore the enclosed volume.

In our problem, we cannot insert the function $u$ in the index
form~\eqref{eq:indexform} since its support is not included in $\Sg$. 
Our first objective is extending the validity of inequality
$Q(u,u)\geq 0$ to more general functions.  We need some preliminary
results.

\begin{remark}
If $\Om$ is a bounded set of finite perimeter in a cone $M$, then
$\Om$ has finite perimeter in $\rr^{n+1}$ by \cite{giusti}, 
Theorem 2.10, and so, the Divergence Theorem is valid in $\Om$ for
vector fields $X\in C^\infty_0(\rr^{n+1},\rr^{n+1})$.  In particular,
if $X(p)\in T_{p}(\ptl M)$ for any $p\in\ptl M-\{0\}$, then
\begin{equation}
\label{eq:gaussgreen}
\int_\Om\text{div}\,X\,d\mathcal{H}^{n+1}=-\int_{\ptl\Om\cap
M}\escpr{X,N}\, d\mathcal{H}^n,
\end{equation}
where $N$ is the generalized inner normal to $\Om$, see 
\cite{bvfunctions}, Theorem~5.8.2.
\end{remark}

\begin{lemma}
\label{lem:monotonicity}
Let $\Om$ be a bounded set of finite perimeter in $M$ with
$\Lambda=\Sg\cup\Sg_0$, where $\Sg_0$ is a closed singular set such
that $\mathcal{H}^{n-2}(\Sg_0)=0$ or consists of isolated points, and
$\Sg$ is a smooth, embedded hypersurface with
$\ptl\Sg-\Sg_{0}\subset\ptl M$.

If $\Om$ is stationary, then there exists a positive constant $\mu$
\emph{(}not depending on $r$\emph{)} such that
\[
\frac{\mathcal{H}^n(\Sg\cap B_r)}{r^n}\leq\mu,
\]
for almost all open balls $B_r\subset\rr^{n+1}$ with small enough
radius.
\end{lemma} 

\begin{proof}
Let $H$ be the constant mean curvature of $\Sg$ with respect to the
normal $N$ which points into $\Om$.  We follow the proofs of the
monotonicity formulae established in \cite{simon} and
\cite{gruterjost}, so we need to show that $\Sg$ has generalized mean
curvature in the sense of \cite{simon}, 17.1, and 
\cite{gruterjost}, p.~133, (7).  That is, we need to establish
\begin{equation}
\label{eq:gmc}
\int_\Sg\text{div}_\Sg\,
X\,d\mathcal{H}^n=-nH\int_\Sg\escpr{X,N}\,d\mathcal{H}^n,
\end{equation}
for vector fields $X$ with support in a small ball (here $\text{div}_\Sg X$ denotes the divergence on $\Sg$).  Fo\-llowing
\cite{simon} and \cite{gruterjost}, the above formula specialize to
appropriate vector fields and the boundedness of the mean curvature of
$\Sg$ give the monotonicity formulae for interior balls 
\cite{simon}, Theorem 17.6, and for balls meeting $\ptl M-\{0\}$, 
see \cite{gruterjost}, Theorem 3.4.  On the other hand, the vector field $X(p)=p$ is
tangent to $\ptl M$ and so, the proof of L.~Simon shows that the
monotonicity formula for interior balls also holds for small balls
centered at the vertex of the cone.  From these formulae the proof of
the Lemma follows.

So let us prove~\eqref{eq:gmc}.  Consider a small open ball $B$
contained in the interior of the cone or centered at a point of
$\Lambda\cap\ptl M$, and a vector field $X\in C^\infty_0(B,\rr^{n+1})$
which vanishes at the origin and such that $X(p)\in T_p(\ptl M)$ for
all $p\in\ptl M-\{0\}$.  Let $\{\varphi_t\}_t$ of $X$ be the one-parameter
group of diffeomorphisms associated to $X$. 

Now, take a point $q_0\in(\Sg\cap M)-\overline{B}$ and a small
Euclidean ball $B'$ disjoint from $B$, and let $W=B'\cap\Sg$.  Let
$u\in C^\infty_0(W)$ such that $-\int_W u\,d\mathcal{H}^n=1$.  It is
possible to construct a vector field $Y$ in $B'$ with $Y|_W=uN$ with
associated flow $\{\psi_s\}_{s}$.

Consider the two-parameter variation $\psi_{s}(\varphi_{t}(\Om))$,
and let $\text{vol}(s,t)=\text{vol}(\psi_{s}(\varphi_{t}(\Om)))$, and
$\mathcal{P}(s,t)=\mathcal{P}(\psi_{s}(\varphi_{t}(\Om)))$. By the 
first variation for volume and \eqref{eq:gaussgreen} we have
\[
\text{vol}_{t}(0,0)=-\int_{\Sg} \escpr{X,N}\,d\mathcal{H}^{n}, \qquad
\text{vol}_{s}(0,0)=-\int_{\Sg}u\,d\mathcal{H}^{n}=1.
\]
By the Implicit Function Theorem, there is a function $s(t)$, with 
$s(0)=0$ and $s'(0)=\int_{\Sg} \escpr{X,N}\,d\mathcal{H}^{n}$, such 
that $\text{vol}(s(t),t)=\text{vol}(\Om)$.

On the other hand, the first variation formula of perimeter for
varifolds \cite{simon}, 16.2, imply
\[
\mathcal{P}_{t}(0,0)=\int_\Sg\text{div}_{\Sg}\,X\,d\mathcal{H}^{n},
\qquad
\mathcal{P}_{s}(0,0)=\int_\Sg\text{div}_{\Sg}\,Y\,d\mathcal{H}^{n}=-\int_{\Sg}
nH\,u\, d\mathcal{H}^{n}=nH.
\]

Since $\Om$ is stationary and the variation
$\psi_{s(t)}(\varphi_{t}(\Om))$ preserves the volume of $\Om$ we get
that the derivative of perimeter with respect to this variation
vanishes.  A direct computation using that the supports of $X$ and $Y$
are disjoint shows that this derivative equals
\begin{align*}
0&=\mathcal{P}_{t}(0,0)+s'(0)\,\mathcal{P}_{s}(0,0) \\
&=\int_\Sg\text{div}_{\Sg}\,X\,d\mathcal{H}^{n}+nH\int_\Sg
\escpr{X,N}\,d\mathcal{H}^{n},
\end{align*}
which proves \eqref{eq:gmc}.
\end{proof}

The next Lemma gives an approximation of the function $1$ on $\Sg$ by 
functions with compact support on $\Sg$ with controlled gradient.

\begin{lemma}
\label{lem:phieps}
Let $\Om$ be a bounded set of finite perimeter in $M$ with
$\Lambda=\Sg\cup\Sg_0$, where $\Sg_0$ is a closed singular set such
that $\mathcal{H}^{n-2}(\Sg_0)=0$ or consists of isolated points,
$\Sg$ is a smooth, embedded hypersurface with boundary, and
$\ptl^{*}\Sg=\ptl\Sg-\Sg_{0}$ is contained in $\ptl M$.

If $\Om$ is stationary then, given $\eps >0$, there is a smooth
function $\varphi_\eps:\Sg\to [0,1]$ with compact support, and such
that
\begin{enum}
\item $\mathcal{H}^n(\{\varphi_\eps\neq 1\})\leq\eps,\quad
\mathcal{H}^{n-1}(\{\varphi_\eps|_{\ptl^{*}\Sg}\neq1\})\leq\eps$. 
\item $\int_\Sg |\nabla\varphi_\eps|^2\,d\mathcal{H}^n \leq\eps$.
\end{enum}
\end{lemma}

\begin{proof}
Suppose that the singular set $\Sg_0$ is nonvoid (otherwise it
suffices to define $\varphi_\eps\equiv 1$ for all $\eps>0$).  Assume
first that $n\geq 3$; in this case it is clear that
$\mathcal{H}^{n-2}(\Sigma_0)=0$ and we proceed exactly as in
\cite{zumbrun2}, Lemma 2.4, to define the function $\varphi_\eps$.  On
the other hand, for $n=2$ the set $\Sg_0$ consists in a finite set of
isolated boundary points and the construction of $\varphi_\eps$ can be
obtained as in \cite{morganritore}, Lemma 3.1.  Both constructions use
Lemma \ref{lem:monotonicity}.
\end{proof}

The following result illustrates how to use Lemma~\ref{lem:phieps} in
order to extend classical results to more general situations.  Recall
that if $\ptl^{*}\Sg=\emptyset$ then the corresponding integral is
assumed to be zero.

\begin{lemma}
\label{lem:divtheorem}
Let $\Om$ be a bounded set of finite perimeter in $M$ with
$\Lambda=\Sg\cup\Sg_0$, where $\Sg_0$ is a closed singular set such
that $\mathcal{H}^{n-2}(\Sg_0)=0$ or consists of isolated points,
$\Sg$ is a smooth, embedded hypersurface with boundary, and
$\ptl^{*}\Sg=\ptl\Sg-\Sg_{0}$ is contained in $\ptl M$.  If $\Om$ is
stationary, then
\begin{enum} 
\item \emph{(}Divergence Theorem\emph{)} If $X$ is a smooth, tangent
vector field on $\Sg$ satisfying $|X|^2$, $\emph{div} X\in L^1(\Sg)$
and $\escpr{X,\nu}\in L^1(\ptl^{*}\Sg)$, then
\[
\int_\Sg\emph{div}\,X\,d\mathcal{H}^n=
-\int_{\ptl^{*}\Sg}\escpr{X,\nu}\,d\mathcal{H}^{n-1},
\]
where $\emph{div}X$ denotes the divergence relative to $\Sg$, and
$\nu$ is the inward normal vector to $\ptl^{*}\Sg$ in $\Sg$.  \item
\emph{(}Integration by parts\emph{)} For any smooth bounded function
$u:\Sg\to\rr$ such that $|\nabla u|^2$, $|\Delta u|\in L^1(\Sg)$ and
$\frac{\ptl u}{\ptl\nu}\in L^1(\ptl^{*}\Sg)$, it is satisfied
\[
\int_\Sg (u\Delta u+|\nabla u|^2)\,d\mathcal{H}^n=
-\int_{\ptl^{*}\Sg}u\,\,\frac{\ptl u}{\ptl\nu}\,\,d\mathcal{H}^{n-1}.
\] 
\end{enum}
\end{lemma}

\begin{proof}
Let us prove (i).  Consider a sequence of functions
$\{\varphi_\eps\}_{\eps>0}$ as obtained in Lemma~\ref{lem:phieps}. 
Since the field $\varphi_\eps X$ has compact support in $\Sg$, we can
apply the classical divergence theorem to obtain
\begin{equation}
\label{eq:divclassic}
\int_\Sg\text{div}(\varphi_\eps X)\,d\mathcal{H}^n=
-\int_{\ptl^{*}\Sg}\varphi_\eps\escpr{X,\nu}\,d\mathcal{H}^{n-1}.
\end{equation}

On the other hand, we have $\text{div}(\varphi_\eps
X)=\varphi_\eps\,\text{div}\,X+\escpr{\nabla\varphi_\eps,X}$.  Letting
$\eps\to 0$ in \eqref{eq:divclassic} and using the Cauchy-Schwarz
inequality in $L^2(\Sg)$, the properties of $\varphi_\eps$, and the
hypotheses on $X$, we prove (i).  Statement (ii) follows from (i) by
taking $X=u\,(\nabla u)$.
\end{proof}

\begin{example}
This example illustrates that the Divergence Theorem for
noncompact manifolds with boundary need not hold if we do not
impose $|X|^2\in L^1(\Sg)$. Let $\Sg$ be the square 
$[0,1]\times [0,1]$ where the vertices have been 
suppressed. The vector field $X(p)=p/|p|^2$ satisfies \cite{whitney}
\[
0=\int_\Sg\text{div}\,X\,d\mathcal{H}^{2}
\neq\int_{\ptl\Sg}\escpr{X,\nu}\,d\mathcal{H}^1=\frac{\pi}{2}.
\]    
\end{example}

Now, we can extend the class of test functions for which inequality 
$Q(u,u)\geq 0$ is valid.

\begin{lemma}
\label{lem:testfunctions}
Let $\Om$ be a bounded set of finite perimeter in $M$ with
$\Lambda=\Sg\cup\Sg_0$, where $\Sg_0$ is a closed singular set such
that $\mathcal{H}^{n-2}(\Sg_0)=0$ or consists of isolated points,
$\Sg$ is a smooth, embedded hypersurface with boundary, and
$\ptl^{*}\Sg=\ptl\Sg-\Sg_{0}$ is contained in $\ptl M$.

If $\Om$ is stable and the cone $M$ is convex, then the index form
defined in \eqref{eq:indexform} satisfies $Q(u,u)\geq 0$ for any
smooth bounded function $u$ on $\Sg$ with mean zero and satisfying
$|\nabla u|^2,|\Delta u|\in L^1(\Sg)$ and $\frac{\ptl u}{\ptl\nu}\in
L^1(\ptl^{*}\Sg)$.
\end{lemma}

\begin{proof} 
Let $u$ be any function under the hypotheses of the statement.  By
Lemma~\ref{lem:divtheorem} we may integrate by parts to conclude that
the index form $Q(u,u)$ in \eqref{eq:indexform} equals $I(u,u)$ in
\eqref{eq:indexform2}.  Hence, the proof finishes if we prove that
$I(u,u)\geq 0$.  To prove this inequality, we can proceed as in
\cite{morganritore}, Lemma 3.3.  We use the sequence
$\{\varphi_\eps\}$ of Lemma~\ref{lem:phieps} to construct a sequence
of functions $\{u_\eps\}_{\eps>0}$ with mean zero and compact support
in $\Sg$, such that $\{u_\eps\}\to u$ in the Sobolev space $H^1(\Sg)$. 
By the stability of $\Om$, we get $I(u_\eps,u_\eps)\geq 0$, which is
equivalent to
\[
\int_\Sg|\sigma|^2\,u_\eps^2\,d\mathcal{H}^n+
\int_{\ptl^{*}\Sg}\text{II}(N,N)\,u_\eps^2\,d\mathcal{H}^{n-1}
\leq\int_\Sg|\nabla u_\eps|^2\,d\mathcal{H}^n.
\]  
Note that $\text{II}(N,N)\geq 0$ since $M$ is convex and so, we can
take $\liminf$ in the above inequality and use Fatou's~Lemma to
deduce
\[
\int_\Sg|\sigma|^2\,u^2\,d\mathcal{H}^n+
\int_{\ptl^{*}\Sg}\text{II}(N,N)\,u^2\,d\mathcal{H}^{n-1}
\leq\int_\Sg|\nabla u|^2\,d\mathcal{H}^n<\infty,
\]
from which the proof follows.
\end{proof}

\begin{remark}
Inequality $I(u,u)\geq 0$ for functions $u$ as in the above lemma was
obtained by P.~Sternberg and K.~Zumbrun~\cite{zumbrun1} for stable
regions with a singular set $\Sg_0$ such that
$\mathcal{H}^{n-2}(\Sg_0)=0$.  However, they did not treat the case of
isolated singularities.
\end{remark}

The following lemma will be useful in proving that the test function
$u=1+H\escpr{X,N}$ satisfies the hypotheses of
Lemma~\ref{lem:testfunctions}.

\begin{lemma}
\label{lem:integrability}
Let $\Om$ be a bounded set of finite perimeter in $M$ with
$\Lambda=\Sg\cup\Sg_0$, where $\Sg_0$ is a closed singular set such
that $\mathcal{H}^{n-2}(\Sg_0)=0$ or consists of isolated points,
$\Sg$ is a smooth, embedded hypersurface with boundary, and
$\ptl^{*}\Sg=\ptl\Sg-\Sg_{0}$ is contained in $\ptl M$.

Let $N$ be the normal along $\Sg$ pointing into $\Om$, $|\sigma|^2$
the squared sum of the principal curvatures associated to $N$ and\/
$\text{\emph{II}}$ the second fundamental form of $\ptl M-\{0\}$ with
respect to the inner normal.

If $\Om$ is stable and $M$ is convex, then $|\sigma|^2\in L^1(\Sg)$
and $\text{\emph{II}}(N,N)\in L^1(\ptl^{*}\Sg)$.
\end{lemma}

\begin{proof}
We reason as in \cite{morganritore}, Lemma 3.3.  Clearly, it suffices
to show that $\int_W|\sg|^2\,d\mathcal{H}^n$,
$\int_{W\cap\ptl^{*}\Sg}\text{II}(N,N)\,d\mathcal{H}^{n-1}<\infty$,
where $W$ is a small open neighborhood of $\Sg_0$ in $\Sg$.  To prove
this we take the function $f\equiv 1$ in $W$ and extend it on $\Sg$ to
a bounded function with mean zero and such that $|\nabla f|^2\in
L^1(\Sg)$.  By applying Lemma~\ref{lem:testfunctions} to $f$, we have
\[
\int_W|\sigma|^2\,d\mathcal{H}^n
+\int_{W\cap\ptl^{*}\Sg}\text{II}(N,N)\,d\mathcal{H}^{n-1}
\leq\int_\Sg|\nabla f|^2\,d\mathcal{H}^n<\infty.
\]
\end{proof}

In order to compute the index form $Q(u,u)$ given by
\eqref{eq:indexform} for $u=1+H\escpr{X,N}$ we need the following
result

\begin{lemma}
\label{lem:partialg}
Let $\Sg$ be a smooth, embedded hypersurface with boundary in a cone
$M$.  Denote by $\ptl^{*}\Sg$ the set of points $x\in\ptl\Sg\cap (\ptl
M-\{0\})$ such that $\Sg$ meets $\ptl M$ orthogonally in a
neighborhood of $x$.  Let $N$ be a unit normal vector to $\Sg$, $\nu$
the inner normal to $\ptl^{*}\Sg$ in $\Sg$, $X(p)=p$ the position
vector field on $\rr^{n+1}$, and $g=\escpr{X,N}$ the support function
of $\Sg$.  Then we have
\[
\frac{\ptl g}{\ptl\nu}=-\text{\emph{II}}(N,N)\,g,\qquad\text{in}\ 
\ptl^{*}\Sg,
\] 
where $\text{\emph{II}}$ is the second fundamental form of $\ptl
M-\{0\}$ with respect to the inner normal.
\end{lemma}

\begin{proof}
Denote by $\nu^*$ the inward normal vector to $\ptl M-\{0\}$.  First
note that $\escpr{X,\nu^*}=0$ in $\ptl M$.  As $\Sg$ meets $\ptl M$
orthogonally in $\ptl^{*}\Sg$, we can differentiate the previous
equality with respect to $N$.  Denoting by $D$ the covariant
derivative in $\rr^{n+1}$ and taking into account that $D_{e}X=e$ for
any vector $e$, we obtain
\begin{equation}
\label{eq:one}
\escpr{X,D_N\nu^*}=0\quad\text{in}\quad\ptl^{*}\Sg.
\end{equation}

Second, we differentiate the equality $\escpr{\nu,N}=0$ with respect
to a vector $v\in T(\ptl^{*}\Sg)$, to~get
\[
\text{II}(N,v)+\text{II}_\Sg(\nu,v)=0,
\]
where $\text{II}_\Sg$ is the second fundamental form of $\Sg$ with
respect to $N$.  We have used the fact that $\nu=\nu^*$ in
$\ptl^{*}\Sg$.  Taking into account the symmetry of $\text{II}$ and
$\text{II}_\Sg$, the above equality becomes
\[
\escpr{v,D_N\nu^*+D_\nu N}=0,\quad\text{for any}\quad v\in T(\ptl^{*}\Sg),
\]
from which we conclude
\[
D_N\nu^*+D_\nu N=a\nu-\text{II}(N,N) N \quad\text{in}\quad
\ptl^{*}\Sg.
\]

Finally, by using \eqref{eq:one} and equality $\escpr{X,\nu}=0$ in
$\ptl^{*}\Sg$, we conclude that
\[
\frac{\ptl g}{\ptl\nu}=D_\nu g=\escpr{X,D_\nu N}=
\escpr{X,a\nu-\text{II}(N,N) N-D_N\nu^*}=-\text{II}(N,N) g.
\]
\end{proof}

Now we can prove the main result of this section

\begin{theorem}
\label{th:stability}
Let $\Om$ be a bounded set of finite perimeter in a convex cone $M$
with $\overline{\ptl\Om\cap M} =\Sg\cup\Sg_0$, where $\Sg_0$ is a
closed singular set such that $\mathcal{H}^{n-2}(\Sg_0)=0$ or consists
of isolated points, $\Sg$ is a smooth, embedded hypersurface with
boundary, and $\ptl^{*}\Sg=\ptl\Sg-\Sg_{0}$ is contained in $\ptl M$.

If $\Om$ is stable, then either $\Om$ is a round open ball in $M$, or
the intersection of a ball centered at the vertex of the cone
with $M$, or a half-ball in $M$ with boundary in a flat piece of $\ptl M$.
\end{theorem}

\begin{proof}
Let $X(p)=p$ be the position vector field on $\rr^{n+1}$, $N$ the
normal to $\Sg$ pointing into $\Om$, and $g=\escpr{X,N}$ the support
function of $\Sg$ with respect to the origin.  Consider the smooth,
bounded function $u$ on $\Sg$, given by
\[
u=1+Hg,
\]
where $H$ is the constant mean curvature of $\Sg$ with respect to $N$.
 
Let us prove that $u$ satisfies the hypotheses of Lemma
\ref{lem:testfunctions}. First note that
\begin{equation}
\label{eq:calculations}
\nabla u=-H\,\sum_i k_i\escpr{X,e_i}e_i, 
\quad\Delta u=|\sigma|^2-nH^2-|\sigma|^2u,
\quad \frac{\ptl u}{\ptl\nu}=-\text{II}(N,N)\,Hg,
\end{equation}
where $e_i$ is a principal direction of $\Sg$ with principal curvature
$k_i$.  The Laplacian $\Delta u$ has been computed from 
\cite{bdc}, Lemma 3.5.  The derivative $\partial u/\partial\nu$ is calculated by
using Lemma \ref{lem:partialg}.  By equation~\eqref{eq:calculations}
and Lemma~\ref{lem:integrability} we obtain $|\nabla u|^2,|\Delta
u|\in L^1(\Sg)$ and $\partial u/\partial\nu\in L^1(\ptl^{*}\Sg )$.

Second, note that $u=n^{-1}\text{div}X^T$, where $X^T(p)$ denotes the
tangent projection of $X(p)=p$ in $T_p\Sg$.  As $\Sg$ is bounded, we
can apply the Divergence Theorem in Lemma~\ref{lem:divtheorem} (i) to
obtain
\begin{equation}
\label{eq:minkowski1}
\int_\Sg u\,d\mathcal{H}^n=-
n^{-1}\int_{\ptl^{*}\Sg}\escpr{X,\nu}\,d\mathcal{H}^{n-1}=0,
\end{equation}
since $X$ is tangent to $\ptl M-\{0\}$ and $\nu$ coincides with the
inner normal vector to $\ptl M-\{0\}$.  Formula \eqref{eq:minkowski1}
is known as the first Minkowski formula.  This formula implies that
$H\neq 0$ and enables us to show that $\Sg$ is connected.  Otherwise,
let $\Sg_1$ and $\Sg_2$ be two components of $\Sg$.  Consider a
locally constant nowhere vanishing function $v$ with mean zero on
$\Sg_{1}\cup\Sg_{2}$.  It is straightforward to check that we can
insert the function $v$ in the index form $Q$ by
Lemma~\ref{lem:testfunctions}.  By the stability of $\Sg$ we get
$Q(v,v)\ge 0$.  But
\[
Q(v,v)=\sum_{i=1,2}\bigg\{-\int_{\Sg_{i}}|\sg|^{2}c_{i}^{2}
-\int_{\ptl^{*}\Sg_{i}} \text{II}(N,N)\,c_{i}^{2}\bigg\}<0,
\]
for some constants $c_{1}$, $c_{2}\neq 0$, a contradiction.

As $u$ is a mean zero function and satisfies the hypotheses of
Lemma~\ref{lem:testfunctions} we can assert that $Q(u,u)\geq 0$.  Now
we compute $Q(u,u)$.  By \eqref{eq:calculations} we have $\Delta
u+|\sg|^2u=|\sg|^2-nH^2$, and so
\begin{equation}
\label{eq:s1}
-\int_\Sg u\,(\Delta u+|\sg|^2\,u)\,d\mathcal{H}^n= -\int_\Sg
(|\sg|^2-nH^2)\,d\mathcal{H}^n-H\int_\Sg
(|\sg|^2-nH^2)\,g\,d\mathcal{H}^n.
\end{equation}

In order to compute the last integral in the above equality 
we use Lemma \ref{lem:divtheorem} (i) with the tangent  
vector field $Y=\nabla g$, and we get
\[
\int_\Sg|\sg|^2\,g\,d\mathcal{H}^n= -\int_\Sg
nH\,d\mathcal{H}^n-\int_{\ptl^{*}\Sg}
\text{II}(N,N)\,g\,d\mathcal{H}^{n-1},
\]
where we have employed Lemma~\ref{lem:partialg} to compute $\frac{\ptl
g}{\ptl\nu}$ and \cite{bdc}, Lemma~3.5, to compute $\Delta g$.  We
deduce, taking into account \eqref{eq:minkowski1}, that
\begin{equation}
\label{eq:s2}
\int_\Sg (|\sg|^2-nH^2)\,g\,d\mathcal{H}^n=-\int_{\ptl^{*}\Sg}
\text{II}(N,N)\,g\,d\mathcal{H}^{n-1},
\end{equation}
which is known as the second Minkowski formula.
On the other hand, the third equality in \eqref{eq:calculations} gives us
\begin{equation}
\label{eq:s3}
\frac{\ptl u}{\ptl\nu}+\text{II}(N,N)\,u=\text{II}(N,N)\quad\text{ in
}(\ptl^{*}\Sg)-\{0\}.
\end{equation}

Finally, by \eqref{eq:s1}, \eqref{eq:s2} and \eqref{eq:s3} we 
conclude
\[
\label{eq:s4}
0\leq Q(u,u)=-\int_\Sg (|\sg|^2-nH^2)\,d\mathcal{H}^n-
\int_{\ptl^{*}\Sg} \text{II}(N,N)\,d\mathcal{H}^{n-1}.
\]

The last inequality implies, by the convexity of the cone and
\cite{bdc}, Lemma 3.2, that
\[
|\sg|^2\equiv nH^2 \ \ \text{and}\quad \text{II}(N,N)\equiv 0,
\]
and so $\Sg$ is a connected piece of a totally umbilical sphere. 
Hence the normal vector $N$ to $\Sg$ can be extended continuously to
all the boundary of $\Om$ and, by \cite{giusti}, Theorem 4.11, we
deduce that $\ptl\Om\cap M$ is a smooth hypersurface.  From a
connectedness argument we conclude that the singular set $\Sg_{0}$ is
empty and, therefore, $\Sg$ is a compact connected piece of a round
sphere. If $\ptl\Sg=\emptyset$ then $\Sg$ is a sphere contained in 
$\overline{M}$. If $\ptl\Sg\neq\emptyset$ the conclusion follows from 
Lemma~\ref{lem:spheres} below. 
\end{proof}

\begin{lemma}
\label{lem:spheres}
Let $M\subset\rr^{n+1}$ be a convex cone and $S\subset\overline{M}$ a
compact, connected subset of a sphere such that $\ptl S\subset\ptl M$
and the boundary of $S$ meets $\ptl M$ orthogonally.  Then, either $S$
is the intersection of a sphere centered at the vertex with the cone
or a half-sphere centered at $\ptl M$ and lying over a flat piece of
$\ptl M$.
\end{lemma}

\begin{proof}
Let $N$ be the inner normal to $S$.  The linear subspace generated by
a vector $e\in\rr^{n+1}$ will be denoted by $L(e)$.  For any $x\in S$,
the normal line $x+L(N(x))$ to $S$ contains the center $x_0$ of $S$. 
Taking a point $y\in\ptl S$ at maximum distance from the vertex of the
cone and using the orthogonality condition, we deduce that $N(y)$ is
proportional to $y$, and so $x_0\in\ptl M\cup (-\ptl M)$.  Cutting
with a plane passing through $0$, $y$ and containing $N(y)$ we see
that the center of the sphere $x_0\in\ptl M$.

If $x_0=0$ then $S$ is a sphere centered at the vertex intersected
with $M$.  If $x_0\neq 0$ then choose any $x\in\ptl S-\{0\}$.  The
tangent hyperplane $\Pi_x$ to $\ptl M$ at $x$ is a supporting
hyperplane of the convex cone $M$ and contains the straight line
$x+L(N(x))$ since $N(x)$ is tangent to $\ptl M$ by the orthogonality
condition.  So $x_{0},x\in\Pi_x\cap\ptl M$ and we conclude that the
segment line $[x_0,x]$ is contained in $\Pi_x\cap\ptl M$ by the
convexity of $M$.  From here we get that $\ptl S\subset\Pi_0\cap\ptl
M$, where $\Pi_0$ is the tangent hyperplane to $\ptl M$ at $x_0$. 
Hence $\ptl S$ is a great circle of $S$, and $\ptl S$ bounds a flat
region in $\Pi_{0}\cap\ptl M$.
\end{proof}

As a consequence of the classification of stable regions in a convex 
cone $M$, we can show which are the isoperimetric regions in $M$

\begin{theorem}[\cite{lions-pacella}]
\label{th:main}
Isoperimetric regions in a convex cone $M\subset\rr^{n+1}$ different
from a half-space are balls centered at the vertex intersected with the
cone.
\end{theorem}  

\begin{proof}
By Proposition~\ref{prop:halfvolume} the existence of isoperimetric
regions for any volume is guaranteed.  Let $\Om$ be an isoperimetric
region.  By Proposition~\ref{prop:boundedness} we know that $\Om$ is a
bounded set of finite perimeter.  Moreover, by
Proposition~\ref{prop:regularity}, the singular set $\Sg_{0}$
satisfies $\mathcal{H}^{n-2}(\Sg_{0}-\{0\})=0$ and so
$\mathcal{H}^{n-2}(\Sg_{0})=0$ if $n\ge 3$ or $\Sg_{0}\subset\{0\}$ if
$n=2$.  As $\Om$ is stable we conclude by Theorem~\ref{th:stability}
that $\Om$ is either a ball, a half-ball lying over a flat piece of
$\ptl M$ or the intersection between a ball centered at the vertex and
the cone.  The proof finishes by using that the relative profile of
the balls centered at the vertex intersected with the cone is better
than the profile of the half-space.
\end{proof}

An special case of a convex cone is a half-space of the Euclidean 
space

\begin{corollary}
\label{cor:stablehalf}
Let $\Om$ be a bounded set of finite perimeter in a half-space
$H\subset\rr^{n+1}$ such that $\overline{\ptl\Om\cap H}$ can be
decomposed as a disjoint union $\Sg\cup\Sg_0$, where $\Sg$ is a
smooth, embedded hypersurface with boundary in the hyperplane $\ptl
H$, and $\Sg_0$ is a closed set of singularities such that
$\mathcal{H}^{n-2}(\Sg_0)=0$ or consisting of isolated points.

If $\Om$ is stable, then $\Om$ is an open round ball contained in
${H}$ or a half-ball centered at $\ptl H$.
\end{corollary}

The above result was obtained in P.~Sternberg and
K.~Zumbrun~\cite{zumbrun1} when the set $\Om$ is a local minimizer of
perimeter with a volume constraint.  They did not use the fact that
any local minimizer is stable.  Their proof consists in reflecting the
local minimizer with respect to $\ptl H$ and use the classical
isoperimetric inequality in $\rr^{n+1}$.  We do not know a previous
proof of Corollary~\ref{cor:stablehalf} in the literature.  Of course
the classification of stable regions in $H$ shows which are the
isoperimetric regions in $H$

\begin{corollary}
\label{cor:half-space2}
Isoperimetric regions in a half-space $H\subset\rr^{n+1}$ are the
half-balls centered at the boundary of the half-space.
\end{corollary}

The arguments we have used in this section are also valid in the
Euclidean space $\rr^{n+1}$.  They are even simpler since boundary
terms do not appear, and allow us to prove 

\begin{corollary}[\cite{zumbrun1}]
\label{cor:euclidean}
Let $\Om$ be a bounded set of finite perimeter in $\rr^{n+1}$ such
that $\ptl\Om=\Sg\cup\Sg_{0}$, where $\Sg$ is a smooth, embedded
hypersurface, and $\Sg_0$ is a closed set of singularities such that
$\mathcal{H}^{n-2}(\Sg_0)=0$.

If $\Om$ is stable, then $\Om$ is a round ball in $\rr^{n+1}$.
\end{corollary}

\begin{remark}
The above result was proved in \cite{bdc} for smooth, compact,
immersed hypersurfaces without singularities.  P.~Sternberg and
K.~Zumbrun proved in \cite{zumbrun1} that
Corollary~\ref{cor:euclidean} holds.
\end{remark}

Since the existence of isoperimetric regions in $\rr^{n+1}$ can be
obtained from general results \cite{morgan}, Theorem~13.4, without
appealing to the isoperimetric inequality in Euclidean space, we get

\begin{corollary}
\label{cor:euclidean2}
Isoperimetric regions in $\rr^{n+1}$ are round balls.
\end{corollary}

\begin{problem}
Consider the Clifford torus $T\subset\mathbb{S}^3$
\[
T=\{(x,y,z,t)\in\rr^4:x^2+y^2=z^2+t^2=1/2\}.
\]
$\mathbb{S}^3-T$ is the union of two domains $C$ and $D$, which are
isometric via the antipodal map and satisfy
$\mathcal{H}^3(C)=\mathcal{H}^3(D)=\mathbf{c}_3/2$, where
$\mathbf{c}_3=\mathcal{H}^3(\mathbb{S}^3)$.  Hence, if $M$ is the cone
over $C$ we know by Proposition~\ref{prop:halfvolume} that there exist
isoperimetric regions in $M$ for any volume.  We cannot apply our
results since $M$ is nonconvex (in fact, at any point of $\ptl
M-\{0\}$ there always are two principal curvatures with opposite
values).
\end{problem}

\begin{example}
Let $M\subset\rr^{n+1}$ be a cone over a smooth domain $C$ in
$\mathbb{S}^n$ with $\mathcal{H}^n(C)\leq\mathbf{c}_n/2$.  We know
from Proposition~\ref{prop:halfvolume} that isoperimetric regions
exist in $M$ for any given volume.  In \cite{lions-pacella}, Remark
1.3, an example of a nonconvex cone in which the balls
centered at the vertex intersected with the cone do not minimize
perimeter with a volume constraint is given.
\end{example}

\section{Appendix: a direct proof of the characterization of
isoperimetric regions in a convex cone}
\label{sec:appendix}
\setcounter{equation}{0}

Using the special form of the isoperimetric profile in a cone, we can
give a direct proof of the characterization of isoperimetric regions
in a convex cone.

\begin{proof}[Proof of Theorem~\ref{th:main}]
Let $\Om$ be an isoperimetric region with volume  $V_{0}>0$ (existence
of $\Om$ is guaranteed by Proposition~\ref{prop:halfvolume}).  The set
$\Om$ is bounded by Proposition~\ref{prop:boundedness} and stationary. 
Denote by $\Sg$ the regular part of $\overline{\ptl\Om\cap
M}$.  Let $N$ be the normal along $\Sg$ which points into $\Om$, and
$H$ the constant mean curvature of $\Sg$ with respect to $N$.  From
statement (iii) in Proposition~\ref{prop:regularity} we know that
the singular set $\Sg_0$ consists of isolated points or satisfies
$\mathcal{H}^{n-2}(\Sg_0)=0$.  Hence, we can consider the functions
$\varphi_{\eps}:\Sg\to\rr$ given by Lemma~\ref{lem:phieps}.

Fix $\eps>0$ and take a vector field $X$ over the cone, such that
$X(p)\in T_p(\ptl M)$ whenever $p\in\ptl M-\{0\}$ and
$\left.X\right|_{\Sg}=\varphi_{\eps}N$.  The flow of diffeomorphisms
$\{f_t\}_t$ of $X$ induces a variation $\Om_t=f_t(\Om)$ of $\Om$. 
Call $V(t)=\vol(\Om_t)$ and $\pp(t)=\pp(\Om_t)$.  By the first
variation formulae for perimeter and volume
\begin{equation}
\label{eq:da/dt}
\pp'(t)=\left.\frac{d\pp}{dt}\right|_{t}=
-\int_{\Sg_{t}} nH_{t}\escpr{X,N_{t}}\,d\mathcal{H}^{n}
-\int_{\Sg_{t}\cap\ptl M}\escpr{X,\nu_{t}}\,d\mathcal{H}^{n-1},
\end{equation}
where $\nu_{t}$ is the inward normal vector to $\Sg_t\cap\ptl M$ in
$\Sg_t=f_t(\Sg)$.  We also obtain by \eqref{eq:gaussgreen}
\begin{equation}
\label{eq:dv/dt}
V'(t)=\left.\frac{dV}{dt}\right|_{t}=-\int_{\Sg_{t}}
\escpr{X,N_t}\,d\mathcal{H}^n,
\end{equation}
so that $V'(0)=-\int_{\Sg}\varphi_{\eps}\,d\mathcal{H}^n<0$.  Hence we
can write $t$ as a function of the volume $V=V(t)$ for $V$ close to
$V_0$.  Let $\pp(V)=\pp[t(V)]$.  By \eqref{eq:da/dt} and
\eqref{eq:dv/dt} we have
\begin{equation}
\label{eq:da/dv}
\frac{d\pp}{dV}\bigg|_{V=V_0}
=\bigg(\int_{\Sg}\varphi_{\eps}\,d\mathcal{H}^n\bigg)^{-1}
\bigg(\int_{\Sg}nH\,\varphi_{\eps}\, d\mathcal{H}^n\bigg)=nH.
\end{equation}

Note that, by the definition of the isoperimetric profile,
$\pp^{(n+1)/n}(V)\geq I_M^{(n+1)/n}(V)$ is sa\-tisfied for $V$ close to
$V_0$; moreover, as $\Om$ is an isoperimetric region, both functions
coincide in $V_0$.  Thus
\begin{equation}
\label{eq:profiles}
\left.\frac{d\pp^{(n+1)/n}}{dV}\right|_{V=V_0}>0,\qquad 
\left.\frac{d^2\pp^{(n+1)/n}}{dV^2}\right|_{V=V_0}\geq 0, 
\end{equation}
where we have used that $I_M^{(n+1)/n}(V)$ is a linear function of $V$
(Proposition \ref{prop:profile}).  By the first equality above and
\eqref{eq:da/dv} we deduce that $H>0$.

Now, we shall compute $d^2\pp^{(n+1)/n}/dV^2$ at $V=V_0$.  A
straightforward calculation gives us
\begin{equation}
\label{eq:second}
\left.\frac{d^2\pp^{(n+1)/n}}{dV^2}\right|_{V=V_0}= \bigg
(\frac{n+1}{n}\bigg)\,\pp(\Om)^{1/n}\left\{\frac{1}{n}\,
\pp(\Om)^{-1}\,\pp'(V_0)^2+\pp''(V_0)\right\},
\end{equation}
so we only have to compute $\pp''(V_0)$.  This calculation requires
second variation of perimeter and volume; a detailed development can
be found in \cite{zumbrun2}, Theorem 2.5.  It is obtained
\begin{align}
\label{eq:d2a/dv2}
\pp''(V_0)=\left.\frac{d^2\pp}{dV^2}\right|_{V=V_0}&=
\bigg (\int_\Sg\varphi_\eps\,d\mathcal{H}^n\bigg )^{-2}  
\\
&\times\left\{\int_\Sg (|\nabla\varphi_\eps|^2
-|\sigma|^2\,\varphi_\eps^2)\,d\mathcal{H}^n- \int_{\Sg\cap\ptl
M}\text{II}(N,N)\,\varphi_\eps^2\,
d\mathcal{H}^{n-1}\right\}\nonumber,
\end{align}
where $|\sigma|^2$ is the squared sum of the principal curvatures
associated to $N$, and $\text{II}$ the second fundamental form of
$\ptl M-\{0\}$ with respect to the inner normal (recall that $N$ is
tangent to $\ptl M$ at the points in $\Sg\cap\ptl M$ since $\Sg$ meets
$\ptl M$ orthogonally).

Taking into account \eqref{eq:da/dv}, \eqref{eq:d2a/dv2},
\eqref{eq:second}, and the second inequality in \eqref{eq:profiles}, 
we get
\begin{align*}
\biggl(\frac{n+1}{n}\biggr)\,\pp(\Om)^{1/n}&\times
\biggl\{\pp(\Om)^{-2}\bigg(\int_{\Sg} nH^2\,d\mathcal{H}^n\bigg)
+\bigg(\int_{\Sg} \varphi_{\eps}\biggr)^{-2}
\\
\notag
&\times\biggl(\int_{\Sg}\big(|\nabla\varphi_{\eps}|^2-|\sg|^2
\,\varphi_{\eps}^2\big)\,d\mathcal{H}^n
-\int_{\Sg\cap\ptl M}
\text{II}(N,N)\,\varphi_{\eps}^2\,d\mathcal{H}^{n-1}\bigg)
\biggr\}\geq 0.
\end{align*}

Now, passing the negative terms of the above equation to the right
of the inequality and taking $\liminf$ when $\eps\to 0$ we obtain,
by Fatou's~Lemma and Lemma~\ref{lem:phieps}
\[
\int_\Sg|\sg|^2\,d\mathcal{H}^n+\int_{\Sg\cap\ptl M}
\text{II}(N,N)\,d\mathcal{H}^{n-1}\leq\int_\Sg
nH^2\,d\mathcal{H}^n<\infty,
\]
from which we conclude that
$\int_\Sg|\sg|^2\,d\mathcal{H}^n$, $\int_{\Sg\cap\ptl M}
\text{II}(N,N)\,d\mathcal{H}^{n-1}<\infty$ and
\[
-\int_{\Sg} (|\sg|^2-nH^2)\,d\mathcal{H}^{n} -\int_{\Sg\cap\ptl M}
\text{II}(N,N)\,d\mathcal{H}^{n-1}\ge 0.
\]

Finally, from the convexity of the cone and \cite{bdc}, Lemma
3.2, we obtain
\[
|\sg|^2\equiv nH^2, \qquad \text{II}(N,N)\equiv 0,
\]
and so $\Sg$ is a union of pieces of totally umbilical spheres of the
same radius $1/H$.  To prove that $\Sg$ is compact and connected we
can proceed as in Theorem~\ref{th:stability}.  Finally, from
Lemma~\ref{lem:spheres} we conclude that $\Om$ is either a ball
contained in $M$, or a half-ball lying over a flat piece of $\ptl M$,
or a ball centered at the origin intersected with the cone.  A simple
comparison of the perimeters of the candidates finally shows that the
last ones are the only isoperimetric regions.
\end{proof}
  
\begin{remark}
The proof of Theorem~\ref{th:main} given in this section is also valid
if $C=\sph^n$, i.e., for the Euclidean space $\rr^{n+1}$.  In this
case all the boundary terms do not appear and the computations are
even simpler.  In this way we get another proof of the isoperimetric
property of balls in Euclidean space.
\end{remark}


\providecommand{\bysame}{\leavevmode\hbox to3em{\hrulefill}\thinspace}
\providecommand{\MR}{\relax\ifhmode\unskip\space\fi MR }
\providecommand{\MRhref}[2]{%
  \href{http://www.ams.org/mathscinet-getitem?mr=#1}{#2}
}
\providecommand{\href}[2]{#2}

\end{document}